\documentclass[10pt, 14paper,reqno]{amsart}
\usepackage{xcolor}

\usepackage{amsmath,amsfonts,amssymb}
\usepackage[breaklinks]{hyperref}
\usepackage{graphicx}
\usepackage{longtable}
\makeatletter
\@namedef{subjclassname@2020}{%
  \textup{2020} Mathematics Subject Classification}
\makeatother

\theoremstyle{plain}
\newtheorem{thm}{Theorem}[section]
\newtheorem{lem}{Lemma}[section]
\newtheorem{cor}{Corollary}[section]
\newtheorem{prop}{Proposition}[section]

\newtheorem{defi}{Definition}[section]

\theoremstyle{proof}

\numberwithin{equation}{section}

\begin{document} 
\title[Tuples with several $n'\text{s}$]{Diophantine triples with the property $D(n)$ for distinct $n$}
\author{Kalyan Chakraborty, Shubham Gupta and Azizul Hoque}
\address{K. Chakraborty @Kerala School of Mathematics, Kozhikode-673571, Kerala, India.}
\email{kalychak@ksom.res.in}
\address{S. Gupta @Harish-Chandra Research Institute,  HBNI, Chhatnag Road, Allahabad-211019, India.}
\email{shubhamgupta@hri.res.in}
\address{A. Hoque @Department of Mathematics, Faculty of Science, Rangapara College, Rangapara, Sonitpur-784505, Assam, India.}
\email{ahoque.ms@gmail.com}
\date{\today}
\maketitle
\begin{abstract} We prove that for every integer $n$, there exist 
 infinitely many $D(n)$-triples which are also $D(t)$-triples for $t\in\mathbb{Z}$ with $n\ne t$. We also prove  that there are infinitely many triples  with the property $D(-1)$ in $\mathbb{Z}[i]$ which are also $D(n)$-triple in $\mathbb{Z}[i]$ for two distinct $n$'s other than $n = -1$ and these triples are not equivalent to any triple with the property $D(1)$. 
  \end{abstract}
\section{introduction}
 A set $\{a_1, a_2, \cdots, a_m\}$ of $m$ distinct non-zero integers is called a Diophantine $m$-tuple with the property $D(n)$ for a given integer $n$, if $a_ia_j + n $ is a perfect square for all $1\leq i< j\leq m$. Note that it is also known as $D(n)$-$m$-tuple or $D(n)$-set and a $D(1)$-set is simply a Diophantine $m$-tuple. The problem of constructing such sets was first studied by Diophantus of Alexandria, who found four rational numbers with this property.  First Diophantine quadruple was found by Fermat, viz. $\{1, 3, 8, 120\}$. Later, Baker and Davenport \cite{BD1969} proved that this quadruple can not be extended to quintuple. In \cite{DU2004}, Dujella proved that there exist at most finitely many quintuples and confirmed that there is no Diophantine sextuple. Recently, He et. al. \cite{HTZ2019} came up with a conclusion that there is no Diophantine  quintuple. Analogous problems for $D(-1)$-$m$-set were investigated in \cite{BCM2020, DF2005,T2015}. 

Let  us look at a general concept of $m$-tuple with the property $D(n)$ in a commutative ring with unity. 
Let  $\mathcal{R}$ be a commutative ring with unity. Given $n \in \mathcal{R}$, 
a set $\{a_1, a_2, \cdots, a_m\} \subset \mathcal{R} \setminus \{0\}$ is called Diophantine $m$-tuple with the property $D(n)$  if $a_ia_j + n $ is a square in $\mathcal{R}$ for all $1 \leq i < j \leq m$.

Let $K$ be an imaginary quadratic field and $\mathcal{O}_K$ its ring of integers. In this setup, there are some interesting results for the existence of Diophantine $m$-tuples with the property $D(n)$ in $\mathcal{O}_K$ with a bound on $m$. In 2019, Ad\v zaga \cite{AD2019} confirmed the non-existence of $m$-tuples with the property $D(1)$ in $\mathcal{O}_K$ for $m \geq 42$. Analogously, Gupta  \cite{GU2021} proved that there is no Diophantine $m$-tuple with the property $D(-1)$ when $m \geq 37$. 

It is interesting to see that $D(n)$-quadruples are related to the representations of $n$ by the binary quadratic form $x^2 -y^2$. In particular, Dujella \cite{DU1993} proved that a $D(n)$-quadruple in $\mathbb{Z}$ exists if and only if $n\not\equiv 2\pmod 4$, that is if and only if $n$ is a difference of two squares of rational integers, up to finitely many exceptions. An analogous result was proved for $\mathbb{Z}[i]$ (see \cite{DU1997}) suggesting that the relation between $D(n)$-quadruples and differences of two squares is not just a coincidence. This motivated Franu\v{s}i\'{c} and Jadrijavi\'{c} \cite{FJ2019} to conjecture that  there exists a quadruple with the property $D(n)$ in a  commutative ring $\mathcal{R}$ with unity  if and only if $n$ can be written as a  difference of two squares in $\mathcal{R}$, except for finitely many $n$. 
This conjecture is proven for $\mathbb{Z}[i]$ (see \cite{DU1997}), rings of integers of $\mathbb{Q}[\sqrt{d}]$ (for certain $d$'s) (see  \cite{FR2008}, \cite{FR2009},  \cite{SO2013}), $\mathbb{Q}[\sqrt[3]{2}]$ (see  \cite{FR2013}) and that of  $\mathbb{Q}[\sqrt{2}, \sqrt{3}]$ (see  \cite{FJ2019}). However, Chakraborty et. al. \cite{CGH-22} shown on the other hand that this conjecture does not hold in $\mathbb{Z}[\sqrt{10}]$.

It is natural to ask whether there exists a Diophantine $m$-tuple which is also a $D(n)$-$m$-tuple for some $n\ne 1$. This question was asked by Kihel and Kihel \cite{KK2001} in 2001 and they conjectured that there is no Diophantine triple which is also $D(n)$-triple for some $n\ne 1$. In 2015, Zhang and Grossman \cite{ZG2015} found a counter example to this claim. They found the triple $\{1, 8, 120\}$ which has the property $D(n)$ for $n = 1, 721$. Motivated by this, several authors worked in this direction. In 2018, Ad\v zaga et. al. \cite{AD2018} proved that there exist infinitely many triples with the property $D(1)$ which has the property $D(n)$ for two distinct $n's$ other than $n = 1$. Subsequently in 2020,  Dujella and Petricevic \cite{DP2020} showed the existence of infinitely many quadruples with the property $D(n)$ for two distinct $n's$, having elements $-1, 7$. They also proved \cite{DUE2020} that there exist infinitely many quadruples with the property $D(0)$ wh
 ich has the property $D(n)$, for two distinct $n'$s other than $n=0$, having all square elements. Recently,  Dujella, Kazalicki and Petri\v cevi\' c \cite{DK2021} proved the existence of infinitely many quintuples with the property $D(0)$ which has also the property $D(n)$, where $n \neq 0$, having all square elements. Before stating our main results, we recall the following definitions. 
\begin{defi}
	A quadruple $\{a, b, c, d\}$ with the property $D(-1)$ in $\mathbb{Z}[i]$ is said to be regular if $d$ can be written as one of the following way
	$$
	d = a + b + c -2abc \pm 2rst,
	$$
	where $r = \sqrt{ab - 1}$, $s = \sqrt{ac - 1}$ and $t = \sqrt{bc - 1}$.
\end{defi}

\begin{defi}
	Let $\mathcal{R}$ be a commutative ring with unity and $n_ 1, n_2 \in \mathcal{R}$. $\mathcal{A} = \{a_1, a_2, \cdots, a_m\}$ and $\mathcal{B} = \{b_1, b_2, \cdots,b_m \}$ are Diophantine $m$-tuple with the property $D(n_1)$ and $D(n_2)$ in $\mathcal{R}$, respectively. We say that $\mathcal{A}$ and $\mathcal{B}$  are equivalent to each other  if there exist a $u \in \mathcal{R}\setminus \{0\}$ such that 
	$$
	\mathcal{A} = \{ub_1, ub_2, \cdots, ub_m\} \text{~~and~~} n_1 = u^2n_2. 
	$$
\end{defi}
We prove the following results:
\begin{thm}\label{thm1}
There exist infinitely many triples with the property $D(n)$ which also has the  property $D(t)$ with $n \neq t$. 
\end{thm}
\subsection*{Remark} An analogous result was proved in \cite[Theorem 4]{AD2018} under the condition that $\{2, a, b, c\}$ is a  quadruple with the property $D(1)$, when $c = 2 + a + b + 4ab \pm 2\sqrt{2a + 1}\sqrt{2b + 1}\sqrt{ab + 1}$. 

In \cite{AD2018}, Ad\v zaga et al. found a family a triples $\{a, b, c\}$ with the property $D(1)$ which has the property $D(n)$ too for two distinct $n$'s, where $n$'s are  $(1/4)(a + b + c)^2 - ab - ac - bc$ and $a + b + c$. If we multiply these triples by $i$ then we get a infinite family of triples with the property $D(-1)$ in $\mathbb{Z}[i]$ which has the property $D(n)$ too for two distinct $n$'s where $n$'s are  $-((1/4)(a + b + c)^2 - ab - ac - bc)$ and $-(a + b + c)$. These triples are equivalent to each other. Here, we see that if a triple $\{a, b, c\}$ has the property $D(-1)$ in $\mathbb{Z}[i]$ as well as two distinct $n$'s, then one value of $n$  can be $(1/4)(a + b + c)^2 - ab - ac - bc$ and the other value of $n$ can be $2 - a - b -c$. This criteria helps to determine a family of triples with the property $D(-1)$ in $\mathbb{Z}[i]$ which are also $D(n)$-triple in $\mathbb{Z}[i]$ for two distinct $n$'s other than $n = -1$. Moreover, these triples are not equivalent to any triple with the property $D(1)$ 

\begin{thm}\label{thm2}
Let $\{2, a, b, c\}$ be a regular quadruple with the property $D(-1)$ in $\mathbb{Z}[i]$. Then the triple $\{a, b, c\}$ has the property $D(n)$ in $\mathbb{Z}[i]$ for $n = n_1, n_2, n_3$ with
\begin{align*}
n_1 &= -1,\\
n_2 &= \dfrac{1}{4}(a + b + c)^2 - ab - ac - bc,\\
n_3 &= 2-a-b-c.
\end{align*}
\end{thm}
\begin{cor}\label{cor1}
There exist infinitely many triples  with the property $D(-1)$ in $\mathbb{Z}[i]$ which have the property $D(n)$, for two distinct $n'$s such that $n \neq -1$ and these triples are not equivalent to any triple having the property $D(1)$.
\end{cor}

\section{Preliminaries}
Let $R$ be the ring of integers of a number field $F$.  We use the standard notation $\square$ to indicate that the expression is a complete square. Assume that $\{a, b, c\}$ is triple  in $R$ with the property $D(n)$ for a given integer $n$. Then
\begin{equation}\label{eq1}
ab + n = \square, ~~~ ac + n = \square, ~~~ bc + n = \square.
\end{equation}
This shows that $x=n$ is an integer solution of the the system of equations
$$ab + x = \square, ~~~ ac + x = \square, ~~~ bc + x = \square.$$
 It is easy to see that the above system of equations has finitely many integer solutions since there are only finitely many integral points on the elliptic curve
\begin{equation}\label{eq2}
E : y^2 = (ab + x)(ac + x)(bc + x).
\end{equation}
The following proposition helps us to reduce the investigation of integer solutions of the above system of equations to the investigation of the induced elliptic curve \eqref{eq2}.
\begin{prop}[\cite{KN1992}, Theorem 4.2]\label{prop1}
Let $K$ be a field with characteristic not equal to 2 and 3. Let
$$
E(K) : y^2 = (x - p)(x - q)(x - r).
$$
Let $P = (x, y) \in E(K)$. Then for some $Q \in E(K)$, $P = 2Q$  if and only if $(x - p), (x - q)$ and  $(x - r)$ all are squares in $K$.
\end{prop}
In order to prove Theorem \ref{thm1} and Theorem \ref{thm2}, we need to find points $P \in E(K)$ satisfying $P = 2Q$ and $x$-coordinate of $P$ lies in $R$ for some $Q \in E(K)$.

\section{Proof of Theorem \ref{thm1}}

Assume that $\{a, b, c\}$ is a triple in $\mathbb{Z}$ with the property $D(n)$. Then  for some $r,s,t\in \mathbb{Z}$, we have
$$
ab + n = r^2, ~~ ac + n = s^2, ~~ bc + n = t^2.
$$
Now the induced elliptic curve is given by 
$$
E(Q): y^2 = (x + ab)(x + ac)(x + bc).
$$
Clearly $P = (0, abc) \in E(\mathbb{Q})$ and thus
$$
2P = \Big(\dfrac{1}{4}(a + b + c)^2 - ab - ac - bc , - abc - \dfrac{1}{8}((a + b + c)^2 - 4ab - 4ac - 4bc)(a + b + c))\Big).
$$
If $a + b + c$ is even, then $x(2P) \in \mathbb{Z}$. Further,  
$$
c = a + b \pm 2\sqrt{ab + n}~~~\text{ if and only if}~~~ \pm \sqrt{ab +n} = \dfrac{c - a - b}{2} ~~~\text{ if and only if}~~~ n = \dfrac{(c - a - b)^2}{4} - ab.
$$  
Hence
$$
c = a + b \pm 2\sqrt{ab + n}~~~\text{ if and only if}~~~ n = x(2P). 
$$
So if we take $c \neq a + b \pm 2\sqrt{ab + n}$, then $n \neq x(2P)$. Thus utilizing Proposition \ref{prop1}), we get, 

\begin{lem}\label{lem1}
A triple $\{a, b, c\}$ with $a+b+c$ even, $c\ne a+b\pm 2 \sqrt{ab+n}$ and with the property $D(n)$ satisfies the property
$D(x(2P))$ such that $n \neq x(2P)$. 
\end{lem}
\subsection{Proof of Theorem \ref{thm1}}
We first consider $n$ even, that is  $n=4m$ or  $4m + 2$ for some $m\in \mathbb{Z}$.  We assume the pair $\{1, x^2-n\}$  with the property $D(n)$, where $x \in \mathbb{Z}$. 


Assume that $\{1, x^2-n\}$ is extendable to a triple $\{1, x^2-n, c\}$ with the property $D(n)$. Then we have 
$$\begin{cases}
c + n = y^2,\\
(x^2 - n)c + n  = z^2,
\end{cases} $$ 
for some $y, z \in \mathbb{Z}$. 
We eliminate $c$ to get 
\begin{equation}\label{eqt1}
z^2 - (x^2 - n)y^2 = n(1- x^2 + n).
\end{equation}
Take $c = 1 + (x^2 - n) + 2x$. Then  $c + n = y^2 = (1 + x)^2$ and $(x^2 - n)c + n  = z^2 = (x^2 + x - n)^2$.
Hence, the equation \eqref{eqt1} has a positive solution, say  $(\alpha, \beta) = (z, y)$. Now, let $x^2 - n = \gamma^2$. Then $(x - \gamma)(x + \gamma) = n$. This imples that $(x - \gamma) = c_1$ and $(x + \gamma) = c_2$, for some $c_1, c_2 \in \mathbb{Z}$ such that $c_1c_2 = n$. Since $c_1, c_2 \leq n$, so $x = \dfrac{c_1 + c_2}{2} \leq n$. So, for a fix $n$, $x^2 - n$ is a square, for only finitely many values of $x$. We can choose $x$ such that $x > n$ for taking $x^2 - n$ square-free. Then the equation \eqref{eqt1} has infinitely many integer solutions $(z, y)$ which are given by  
$$
z + \sqrt{x^2 - n}y = (\alpha + (\sqrt{x^2 - n})\beta)(\alpha_1 + (\sqrt{x^2 - n})\beta_1)^t ,~~ t> 0,
$$
where $\alpha_1 + (\sqrt{x^2 - n})\beta_1$ has norm $1$ with $\alpha_1, \beta_1>0$.
Now we put one more restriction to $x$ such that $x^2 - n$ is odd. This implies that $c$ is even and thus both $\alpha$ and $\beta$ are even. Hence for varying $t > 0$, we get infinitely many $z$ which are even.


Utilizing these solutions (i.e., $z$), we can generate infinitely many even $c$ varying $t$, except for finitely many $t$.  Thus, there are infinitely many triples of the form $\{1, x^2-n, c\}$ with the property $D(n)$ which also satisfy the hypothesis of Lemma \ref{lem1}. Hence there are infinitely many triples with  the property $D(n)$, which have the property $D(x(2P))$. 

Secondly, we consider $n = 4m + 1$ or $4m + 3$ for some $m\in\mathbb{Z}$. 
As in the previous case, we arrive at \eqref{eqt1}, $(\alpha, \beta) = (z, y)$ is a positive integer solution of \eqref{eqt1} when $c = 1 + (x^2 - n) + 2x$ and, $x^2 - n$ is odd, by choosing appropriate $x$. Thus, $c$ is even and hence both $\alpha$ and $\beta$ are odd. Therefore, there are infinitely many solutions $(z, y)$ of \eqref{eqt1} of the form 
$$
z + (\sqrt{x^2 - n})y = (\alpha + (\sqrt{x^2 - n})\beta)(\alpha_1 +( \sqrt{x^2 - n})\beta_1)^{2t} , ~~~t > 0,
$$
where $\alpha_1 + \sqrt{x^2 - N}\beta_1$ has norm one with $\alpha_1, \beta_1>0$ ($\alpha_1$ is odd and $\beta_1$  is even).  For varing $t > 0$, we get infinitely many $z$ which are odd and these values of $z$ help us to get infinitely many even $c$ varying $t$, except for finitely many $t$. Therefore, by Lemma \ref{lem1} we can conclude that there are infinitely many triple of the form $\{1, x^2-n, c\}$ with the property $D(n)$ as well as $D(x(2P))$. 

\subsection{Toy examples} 
We first consider $n=4$ and the pair $\{1, 5\}$ with  the property $D(4)$. We need to construct $c$ satisfying the hypothesis of Lemma \ref{lem1}. Let $\{1,  5, c\}$ be a triple with the property $D(4)$.  Then we can find $y, z \in \mathbb{Z}$ such that 
$$\begin{cases}
c + 4 = y^2,\\
5c + 4 = z^2.
\end{cases}
$$
Eliminating $c$, we get 
$$
z^2 - 5y^2 = -16.
$$
Clearly, $(z, y)=(8, 4)$ is solution of the above equation, and thus 
the other solutions are given by positive solutions $(z, y)$ are given by 
\begin{equation}\label{eqe1}
z + y\sqrt{5} = (8 + 4\sqrt{5})(9 + 4\sqrt{5})^t, ~~t \geq 1.
\end{equation}
For $t=1$, we get $y=68$ and thus $c=4620$. 

Clearly, $P=(0, 5\times4620)$ is an integral point on the elliptic curve:
$$E(Q): Y^2=(X+5)(X+4620)(X+23100)$$
and thus $x(2P)=\frac{1}{4}\times 4626^2-5-4620-5\times 4620=5322244$. 

Now $$\begin{cases}
5+5322244=2307^2,\\
4620+5322244=2308^2,\\
5\times4620+5322244=2312^2.
\end{cases}$$
Thus, $\{1, 5, 4620\}$ is both $D(4)$ as well as $D(5322244)$ sets. Similarly each value of $t>0$, we can find a triple of $D(n)$ set for two distinct values of $n$.  

%

Secondly we consider $n = 3$ and the pair $\{1, 13\}$. By using Lemma \ref{lem1}, we determine $c$ such that $\{1, 5, c\}$ makes a triple with the property $D(3)$. So
$$
c + 3 = y^2 \text{~~and~~} 13c + 3 = z^2,
$$
for some $y, z \in \mathbb{Z}$. Eliminating $c$ from above equation, we get
$$
z^2 - 13y^2 = -36.
$$
Above equation has one solution $(z, y) = (17, 5)$. Then we make some positive solution $(z, y)$ such that 
$$
z + \sqrt{13}y = (17 + 5\sqrt{13})(649 + 180\sqrt{13})^{2t}, ~~t \geq 1.
$$
If $t = 1$, then $z = 29507417$, so $c = 66975973693222$. We have $P = (0, 870687658011886)$ and then 
\begin{align*}
x(2P) &= \dfrac{1}{4}(1 + 13 + 66975973693222)^2 - 13 - 66975973693222 - 13 \times 66975973693222\\
 &= 1121445263038322515826332803.
\end{align*}
Now
$$\begin{cases}
1121445263038322515826332803 + 13 \times 1 =33487986846604^2,\\
1121445263038322515826332803 + 66975973693222\times1=33487986846605^2,\\
1121445263038322515826332803 + 66975973693222\times13=33487986846617^2.
\end{cases}$$
Thus a triple $\{1, 13, 66975973693222\}$ has the property $D(3) $ as well as \\ $D(1121445263038322515826332803)$. Similarly each each value of $t>0$, we can find a triple of $D(n)$ set for two distinct values of $n$.
\section{Proof of Theorem \ref{thm2}}
 
Let $\{a, b, c\}$ be a Diophantine triple with the property $D(-1)$ in $\mathbb{Z}[i]$. Then there exist $r, s, t \in \mathbb{Z}[i]$ such that
$$\begin{cases}
ab - 1 = r^2, \\
 ac - 1 = s^2,\\ 
  bc - 1 = t^2.
\end{cases}$$
Assume that the above triple satisfies the property $D(n)$ for $n\in \mathbb{Z}[i]\setminus\{-1\}$. Then 
$$\begin{cases}
ab + n = \square,\\
 ac + n = \square, \\
  bc + n = \square.
\end{cases}
$$
Now we can associate this system of equations with the following elliptic curve over $K=\mathbb{Q}(i)$:
$$
E(K):~y^2 = (x+ab)(x+ac)(x+bc).  
$$
It is easy to see that the points 
$
P = (0, abc)$ and  $S = (-1, \pm rst)\in E(K)\cap (\mathbb{Z}[i]\times \mathbb{Z}[i])$. Thus, if $a+b+c$ is a multiple of $2$ in $\mathbb{Z}[i]$ and $c\ne a+b\pm 2\sqrt{ab+n}$, then by Lemma \ref{lem1} the triple $\{a, b, c\}$ also has the property $D(x(2P))$. 

Assume that $
R = (rs + rt + st - 1, \pm(r + s)(r + t)(s + t))\in E(K)\cap (\mathbb{Z}[i]\times \mathbb{Z}[i])
$. Then $S=2R$ and thus $S + 2P \in 2E(K)$. Now,
\begin{align*}
x(S + 2P) = \dfrac{-1}{4}(a + b + c)^2 + 1 + \dfrac{(8abc + ((a + b + c)^2 - 4ab - 4bc - 4ac)(a + b + c) \pm 8rst)^2}{4((a + b + c)^2 - 4ab - ac - 4bc + 4)^2}.
\end{align*}
By computer search, we have many triples with the property $D(-1)$ in $\mathbb{Z}[i]$ which has $x(S + 2P) = 2- a -b -c$. So let us assume that
$x(S + 2P) = 2 - a - b - c.$ By straightforward calculation, it is equivalent to $e_1e_2e_3 = 0$, where
\begin{align*}
e_1 &= a^2 - 2ab + b^2 - 2ac - 2bc + c^2 + 4,\\
e_2 &= 8abc + a^2 - 2ab + b^2 - 2ac - 2bc + c^2 - 4a - 4b - 4c + 8,\\
e_3 &= a^4 - 2a^2b^2 + b^4 - 2a^2c^2 - 2b^2c^2 + c^4 - 2a^3 + 2a^2b + 2ab^2 - 2b^3 + 2a^2c + 4abc + 2b^2c + 2ac^2  \\ &+ 2bc^2 - 2c^3 + a^2 - 2ab + b^2 - 2ac - 2bc + c^2 - 4a - 4b - 4c + 8.
\end{align*} 
It is easy to see that $e_1 = 0$ if and only if $c = a + b \pm 2\sqrt{ab - 1}$. This does not help us to generate infinitely many triples with the property $D(-1)$ in $\mathbb{Z}[i]
$ which has also the property $D(n)$, for two distinct $n's$ such that $n \neq -1$. 

Let
\begin{equation}\label{eqq1}
e_3 = 0 = f_1^4 - 2f_1^3+ f_1^2 - 4f_1 - 8f_3 - 4f_2 - 4f_2f_1^2 + 8f_1f_2 + 8f_1f_3 + 8,
\end{equation}
where
$f_1 = a + b + c,
f_2 = ab + ac + bc,
f_3 = abc.$ \\
We can simplify $e_3$ in the following way
$$
e_3 = (f_1^2 - 4f_2)(f_1 - 1)^2 + (8f_3 - 4)(f_1 - 1) + 4.
$$
If $e_3 = 0$, so
\begin{equation}\label{eq3}
e_3 = (f_1 - 1)((f_1 - 1)(f_1^2 - 4f_2) + (8f_3 - 4)) = -4.
\end{equation}
After taking norm, we get 
$$
||f_2 - 1|| = 1, 2, 4, 8 \text{~~or~~} 16,
$$
where $||\cdot||$ denote the norm function.
Since we want generate triples with the property $D(-1)$ in $\mathbb{Z}[i]$ with three $n's$, we take $a + b + c$ even, so $||a + b + c -1|| = 1$ (by reducing at modulo 4), and this implies that $a + b + c = 0 \text{~~or~~} 2$.
First consider $a + b + c = 0$, then $c = - a - b$. Put this value of $c$ in the equation \eqref{eqq1}, we get
\begin{equation}\label{eqq2}
2a^2b + 2ab^2 + a^2 + ab + b^2 + 1 = 0.
\end{equation}
This implies that 
$$
b = -\frac{a(2a +  1) \pm \sqrt{4a^4 - 4a^3 -3a^2 - 8a + 4}}{2(2a + 1)}.
$$
Also $2a + 1 \neq 0$, for $a \in \mathbb{Z}[i]$. We can write  $(4a^4 - 4a^3 -3a^2 - 8a + 4) = (2a + 1)^2(a^2 - 2a + 1) - 5(2a + 1)$. Now, $b$ can be written as 
\begin{align*}
b = -\dfrac{a}{2} \pm \dfrac{\sqrt{((a - 1)^2 - 5(2a + 1))/(2a + 1)^2}}{2}
\end{align*} 
So necessary condition for $b \in \mathbb{Z}[i]$ is $\dfrac{5}{(2a + 1)} \in \mathbb{Z}[i]$. So $||2a + 1|| = 1, 5 \text{~~or~~} 25$. This implies that $a = 0, -1, 2, -3$, and then from \eqref{eqq2}, $(a, b) = (-1, 2), (2, -1)$. But these values of $a$ and $b$ does not make pair with the property $D(-1)$  in $\mathbb{Z}[i]$. Hence $a + b + c \neq 0$ so $e_3 = 0$. Now we consider $c = 2 - a - b$. Then equation \eqref{eqq2} will be 
\begin{equation}\label{eqq3}
-2a^2b - 2ab^2 + a^2 + 5ab + b^2 - 2a - 2b = 0.
\end{equation}
This implies that 
$$
b = -\dfrac{(2a - 1)(a - 2) \pm \sqrt{4a^4 - 12a^3 + 13a^2 - 12a + 4}}{2(2a - 1)}.
$$

We have $4a^4 - 12a^3 + 13a^2 - 12a + 4 = (2a - 1)^2(a- 1)^2 -3(2a - 1)$. Similarly as above, For $b \in \mathbb{Z}[i]$, $\dfrac{3}{2a - 1} \in \mathbb{Z}[i]$. This implies that $||2a - 1|| = 1, 3 \text{~~and~~} 9$. From here, we get $a = 1, -1, 2$, and then $(a, b) = (2, 0)$, by using \eqref{eqq3}. From this values of $a$ and $b$ pair $\{a, b\}$  with the property $D(-1)$  in $\mathbb{Z}[i]$ does not make. Therefore $e_3 \neq 0$.


Observe that $e_2 = 0$ if and only if
$$
c = 2 + a+ b - 4ab \pm 2\sqrt{ab - 1}\sqrt{2a - 1}\sqrt{2b - 1}.
$$  
This holds only when $\{2, a, b , c\}$ is a regular quadruple with the property $D(-1)$ in $\mathbb{Z}[i]$. 

Assume that $\{2, a, b , c\}$ is a regular quadruple with the property  $D(-1)$ in $\mathbb{Z}[i]$. Then the triple $\{a, b, c\}$ has the property $D(n)$, for $n =  n_1, n_2$ with
\begin{align*}
n_1 &= x(2P),\\
n_2 &= 2 -a -b -c,
\end{align*}
when the following conditions hold:\\
\begin{itemize}
\item[(i)] $2 + a+ b - 4ab \pm 2\sqrt{ab - 1}\sqrt{2a - 1}\sqrt{2b - 1} \neq a + b \pm 2\sqrt{ab - 1}$,\\
\item[(ii)] $n_1 \neq n_2,$\\
\item[(iii)] $2 - a - b - c \neq -1$,\\
\item[(iv)] $a + b + c$ is a multiple of $2$ in $\mathbb{Z}[i]$.
\end{itemize}
Here, condition $(i)$ and $(iv)$ shows that triple $\{a, b, c\}$ satisfies the hypothesis of Lemma \ref{lem1}, and the condition  $(ii)$ and $(iii)$ tells that all the three $n$'s must be distinct. Now we check the valadity of all these four conditions. 

Clearly (i) holds; otherwise 
$$
4a^2b^4 - 8(a^3 + a)b^3 + 4(a^4 + 4a^2 + a + 1)b^2 + 4a^2 - 4(2a^3 - a^2 + 3a + 1)b - 4a + 5 = 0,
$$
an impossibility. 

Now if (ii) does not hold then 
$$
4a^2b^4 - 8(a^3 - 2a^2 + 2a)b^3 + 4(a^4 + 4a^3 - 8a^2 + 4)b^2 + 16a^2 - 16(a^3 - 3a + 2)b - 32a + 16 = 0.
$$
Then 
$$
\small a =  \dfrac{\pm (2(b - 1)(-2b + 1)^{3/4} + (b^2 - 2b + 2)(-2b + 1)^{1/4} + \sqrt{\pm (8b^3 - 28b^2 +  28b - 8) + B\sqrt{-2b + 1} }\sqrt{b}}{2b(-2b + 1)^{1/4}}
$$
with $B=(b^3 - 12b^2 + 20b - 8)$. 
Since $\{2, b\}$ is a pair with the property $D(-1)$ in $\mathbb{Z}[i]$, so there exist $y \in \mathbb{Z}[i]$, such that $-2b + 1 = y^2$, or $b = (1 - y^2)/2$. Put this value of $b$ in above expression, we get
\begin{equation}\label{eqq4}
a = -\dfrac{y^3 - 3y^2 - y - 5 \pm \sqrt{y^3 - 5y^2 + 3y - 7} \sqrt {y + 1} (y - 1)}{4(y + 1)}.
\end{equation}
or
\begin{equation}\label{eqq5}
a = -\dfrac{y^3  + 3y^2 - y + 5 \pm \sqrt{y^3 + 5y^2 + 3y + 7} \sqrt {y - 1} (y + 1)}{4(y - 1)}.	
\end{equation}
First, let us take equation \eqref{eqq4}. We can write $y^3 - 3y^2 - y -5 = (y^2 - 4y + 3)(y + 1) - 8$. Using this, equation \eqref{eqq4} will be 
$$
(y + 1)(4a + y^2 -4y + 3) + 8 = \pm (\sqrt{y^3 - 5y^2 + 3y - 7} \sqrt {y + 1} (y - 1))   
$$
By squaring, arranging the terms , and after that taking norm both side, we get  that $||y + 1||$ must divide $2^{12}$. From here, we get the values of $b$, then by using equation \eqref{eqq4} we get $a$, and then one can check that $\{a, b\}$ does not make a pair with the property $D(-1)$ in $\mathbb{Z}[i]$.
Now we pick equation $\eqref{eqq5}$. We can write $ = (y^2 + 4y + 3)(y - 1) + 8$, and then follow the above path, and then get same result. Hence, for analying both equation,
$||y + 1||$ must divide $2^{12}.$ From here, we get the values of $b$, then by using equation \eqref{eqq4} we get $a$, and then one can check that $\{a, b\}$ does not make a pair with the property $D(-1)$ in $\mathbb{Z}[i]$. Hence condition $(ii)$ holds.

  
Suppose that
$
a + b + c = 3.
$ Then 
$$
3 = 2(1 +a + b - 2ab \pm \sqrt{ab - 1}\sqrt{2a - 1}\sqrt{2b - 1}),
$$
which is contradiction. Thus (iii) does not hold. 

Finally, $(iv)$ follows from the assumption that $\{2, a, b , c\}$ is a regular quadruple with the property $D(-1)$ in $\mathbb{Z}[i]$. This completes the proof. 

 \section{Proof of the corollary \ref{cor1}}

First we define $a, b$ for a given $m\in \mathbb{Z}[i]$ as follows: 
\begin{align*}
a &= 2m^2 +  2m + 1\\
b &= 2m^2 + 6m + 5,
\end{align*}
so that $\{2, a, b\}$ makes a triple with the property $D(-1)$ in $\mathbb{Z}[i]$. Now we take $c = -32m^4- 128m^3 - 184m^2 - 112m - 24$ so that $\{2, a, b, c\}$ makes a regular quadruple, except for $m=-1$. Then $a, b ,c$ satisfy the assumptions of Theorem \ref{thm2}  and thus the triple $\{a, b, c\}$ has the property $D(n)$ in $\mathbb{Z}[i]$, for $n = n_1, n_2, n_3$ with
\begin{align*}
n_1 &= -1,\\
n_2 &= 32m^4 + 128m^3 + 180m^2 + 104m + 20,\\
n_3 &= 256m^8 + 2048m^7 + 7104m^6 + 13952m^5 + 16992m^4 + 13184m^3 + 6396m^2 + 1784m + 220.
\end{align*}
Now let us discuss another example which is useful to get a triple four distinct $n$'s.\\
\subsection*{An example} By using Theorem \ref{thm2}, a triple $\{2m(m + i), 2(m + 1)(m + 1 + i), -32m^4 + (-64i - 64)m^3 + (-96i + 8)m^2 + (-24i + 40)m + 4i + 8\}$ with the property $D(-1)$ in $\mathbb{Z}[i]$ also has the property $D(n)$, for $n = n_1, n_2$, where
\begin{align*}
n_1 =~~ &(2m + i + 1)^2(8m^3 + 12i + 8)m^2 + (8i - 4)m - 1)\\ &(8m^3 + 12i + 16)m^2 + (16i + 4)m + 4i - 3),\\
n_2 = ~~ &(2m + i + 1)^2(8m^2 + 8(i + 1)m + 4i - 3),
\end{align*}
for $m \neq 0, -i, -1, -1 - i.$
\section{Triple for four distinct $n$'s}
In \cite{AD2018}, it was found that $\{4, 12, 420\}$ is a Diophantine triple with the property $D(n)$ for $n = 1, 436, 3796, 40756$. It is easy to see that $\{4i, 12i, 420i\}$ is a Diophantine triple with the property $D(n)$ in $\mathbb{Z}[i]$ for $n = -1, -436, -3796, -40796$. Cleary, these triples are equivalent to each other.
Now we construct an example of a triple with the property $D(-1)$ which has  the property $D(n)$ for three distinct $n$'s other than $n = -1$ and it is also not equivalent to any triple with the property $D(n)$.
Take $m = -2$ in the example 3. We get the set $\{-4i + 8, -2i + 2, 180i - 40\}$.  Thus, the elliptic curve induced by this triple over $\mathbb{Q}[i]$ is
\begin{align*}
y^2 = ~~&(x + ((-i)(i + 1)^7(2i + 1)))(x + (i(i - 4)(-i - 2)^2(i + 1)^8(2i + 1)^2)) \\ &(x + (-1(i - 4)(-i - 2)^2(i + 1)^7(2i + 1))).
\end{align*}
It is easy to check that $X = (-440i - 120, -2240i + 9280)$ lies on the above curve and thus
\begin{align*}
2X = (-1194i + 392, 24618i - 23374). 
\end{align*}
Hence the triple has the property $D(n)$ for four distinct $n'$s which are as follows:
\begin{align*}
n_1 &=  -1,\\
n_2 &= (-i - 2)^2(i + 1)^2(i + 4)(4i + 9)(11i + 20) = -4626i - 8032,\\
n_3 &= -i(-13i - 12)(-i - 2)^2(i + 1)^2 = -174i + 32,\\
n_4 &= -1194i + 392.
\end{align*}


 

\section*{Acknowledments}
The third author is supported by SERB MATRICS Project (No. MTR/2021/000762), Govt. of India.

\end{document}